\documentclass[10pt,leqno]{article}
\usepackage{amsthm,amsfonts,amssymb,amsmath,oldgerm}
\usepackage{epsfig}
\numberwithin{equation}{section}
\usepackage[thinlines]{easybmat}

\newcommand{\DDD}{D3}

\def\eps{\varepsilon }

\newcommand\R{\mathbb R}

\def\eps{\varepsilon}

\newcommand\errfn{\textrm{errfn}}
\newcommand\br{\begin{remark}}
\newcommand\er{\end{remark}}
\newcommand\bp{\begin{pmatrix}}
\newcommand\ep{\end{pmatrix}}
\newcommand\be{\begin{equation}}
\newcommand\ee{\end{equation}}
\newcommand\ba{\begin{equation}\begin{aligned}}
\newcommand\ea{\end{aligned}\end{equation}}

\newcommand{\bap}{\begin{app}}
\newcommand{\eap}{\end{app}}
\newcommand{\begs}{\begin{exams}}
\newcommand{\eegs}{\end{exams}}
\newcommand{\beg}{\begin{example}}
\newcommand{\eeg}{\end{exaplem}}
\newcommand{\bpr}{\begin{proposition}}
\newcommand{\epr}{\end{proposition}}
\newcommand{\bt}{\begin{theorem}}
\newcommand{\et}{\end{theorem}}
\newcommand{\bc}{\begin{corollary}}
\newcommand{\ec}{\end{corollary}}
\newcommand{\bl}{\begin{lemma}}
\newcommand{\el}{\end{lemma}}
\newcommand{\bd}{\begin{definition}}
\newcommand{\ed}{\end{definition}}
\newcommand{\brs}{\begin{remarks}}
\newcommand{\ers}{\end{remarks}}

\newtheorem{theo}{Theorem}[section]
\newtheorem{prop}[theo]{Proposition}

\newtheorem{exams}[theo]{Examples}

\numberwithin{equation}{section}

\newcommand{\RR}{{\mathbb R}}

\newtheorem{theorem}{Theorem}[section]
\newtheorem{proposition}[theorem]{Proposition}
\newtheorem{corollary}[theorem]{Corollary}
\newtheorem{lemma}[theorem]{Lemma}
\newtheorem{definition}[theorem]{Definition}

\newtheorem{example}[theorem]{Example}
\newtheorem{remark}[theorem]{Remark}

\newcommand{\RM}{\mathbb{R}}

\title{
Stability of periodic Kuramoto--Sivashinsky waves}

\author{\sc \small
Blake Barker\thanks{ \tiny Indiana University, Bloomington, IN 47405;
bhbarker@indiana.edu: Research of B.B. was partially supported
under NSF grants no. DMS-0300487 and DMS-0801745.}
~~~
Mathew A. Johnson\thanks{ \tiny Indiana University, Bloomington, IN 47405;
matjohn@indiana.edu: Research of M.J. was partially supported by an NSF Postdoctoral Fellowship under NSF grant DMS-0902192.}
~~~
Pascal Noble\thanks{ \tiny Universit\'e Lyon I, Villeurbanne, France;
noble@math.univ-lyon1.fr:
Research of P.N. was partially supported by the French ANR Project no.
ANR-09-JCJC-0103-01.}
~~~
L.Miguel Rodrigues\thanks{ \tiny Universit\'e de Lyon, Universit\'e Lyon 1,
Institut Camille Jordan, UMR CNRS 5208, 43 bd du 11 novembre 1918,
F - 69622 Villeurbanne Cedex, France; rodrigues@math.univ-lyon1.fr: Stay of M.R. in Bloomington was supported by 
Frency ANR project no. ANR-09-JCJC-0103-01.}
~~~
Kevin Zumbrun\thanks{ \tiny Indiana University, Bloomington, IN 47405;
kzumbrun@indiana.edu:
Research of K.Z. was partially supported
under NSF grant no. DMS-0300487.}
}
\begin{document}

\maketitle


\begin{center}
\end{center}

\begin{center}
{\bf 2000 MR Subject Classification}: 35B35.
\end{center}


\begin{abstract}
\begin{small}
In this note, we announce a general result resolving the long-standing
question of nonlinear modulational stability, or stability
with respect to localized perturbations, of periodic
traveling-wave solutions of the generalized Kuramoto--Sivashinski
equation, establishing that spectral modulational stability, 
defined in the standard way,
implies nonlinear modulational stability with sharp rates of decay. 
The approach extends readily to other
second- and higher-order parabolic equations, for example, the Cahn Hilliard
equation or more general thin film models.
\end{small}
\end{abstract}


\bigbreak

\section{Introduction }\label{intro}

In this note, we describe recent results obtained
using techniques developed in \cite{JZ,JZN,BJNRZ2} on
linear and nonlinear stability of periodic traveling-wave solutions of the
generalized Kuramoto-Sivashinsky (gKS) equation 
\be\label{e:KS}
u_t+\gamma \partial_x^4u+\epsilon \partial_x^3 u +
\delta \partial_x^2u+\partial_x (u^2/2) =0,
\quad
\delta >0,
\ee
a canonical model for pattern formation in one spatial dimension
that has been used to describe, variously, plasma instabilities,
flame front propagation, turbulence in reaction-diffusion systems,
and thin film flow down an incline \cite{S1,KT,CD,DSS,PSU}.

More generally, we consider (taking without
loss of generality $\gamma=1$) an equation of the form
\be\label{e:gen}
u_t+\partial_x^4u+ \epsilon \partial_x^3u+ \delta \partial_x^2u+\partial_x f(u)=0,
\ee
$f\in C^2$, $\delta$ not necessarily positive.
Our methods apply also, with slight modifications to accomodate
quasilinear form (see \cite{JZN}) to the Cahn Hilliard equation and
other fourth-order models for thin film flow as discussed for
example in \cite{BMSZ}.
Indeed, the argument, and results, extend to arbitrary $2r$-order parabolic
systems, so is essentially completely general for the diffusive case.
As shown in \cite{JZ}, they can apply also to mixed-order and
relaxation type systems in some cases as well.

It has been known since 1976,
almost since the introduction of the model \eqref{e:KS} in 1975 \cite{KT,S1},
that there exist spectrally stable bands of solutions in parameter space; 
see for example the numerical studies in \cite{CKTR,FST}.
Moreover, numerical time-evolution experiments described for
example in \cite{CD} suggest that these waves are nonlinearly
stable as well, serving as attractors in the chaotic dynamics
of (gKS).
However, up to now this conjecture had not been rigorously verified.

Here, we announce the result, resolving this open question, 
that {\it spectral modulational stability},
defined in the standard sense of the modulational 
stability literature,\footnote{
In particular, the sense verified in \cite{CKTR,FST}}
{\it implies linear and nonlinear modulational stability}.
Our analysis gives at the same time new understanding even
at the formal level of Whitham averaged equations.
Further details, along with numerical investigations of existence
and spectral stability, will be given in \cite{BJNRZ3}.

\section{The traveling-wave equation}
Substituting $u=\bar u(x-ct)$ into \eqref{e:gen},
we find that the traveling-wave equation is
$-c u'+u''''+\epsilon u'''+ \delta u'' + f(u)' =0$,
or, integrating once in $x$,
\be\label{e:inttw}
-c u+u'''+\epsilon u''+ \delta u' + f(u) =q,
\ee
where $q\in \R$ is a constant of integration.
Written as a first-order system in $(u,u',u'')$, this is
\be\label{e:1tw}
\bp u\\u'\\u''\ep'=
\bp u'\\
u''\\
c - \epsilon u''- \delta u' - f(u) +q
\ep.
\ee

It follows that periodic solutions of \eqref{e:gen}
correspond to values
$(X,c,q,b)\in \RR^6$, where $X$, $c$, and $q$ denote period,
speed, and constant of integration, and $b=(b_1,b_2,b_3)$ denotes
the values of $(u,u',u'')$ at $x=0$, such that
the values of $(u,u',u'')$ at $x=X$ of the solution of
\eqref{e:1tw} are equal to the initial values $(b_1,b_2,b_3)$.

Following \cite{JZ}, we assume:

(H1) $f\in C^{K+1}$, $K\ge 4$.

(H2) The map $H: \,
\R^6  \rightarrow \R^3$	
taking $(X,c,q,b) \mapsto (u,u',u'')(X,c,b; X)-b$
is full rank at $(\bar{X},\bar c, \bar b)$,
where $(u,u',u'')(\cdot;\cdot)$ is the solution operator of \eqref{e:1tw}.

By the Implicit Function Theorem,
conditions (H1)--(H2) imply that the set of periodic solutions
in the vicinity of $\bar U$ form a
smooth $3$-dimensional manifold (counting translations)
\be\label{manifold}
\{\bar U^\beta(x-\alpha-c(\beta)t)\},
\;
\hbox{\rm with $\alpha\in \RR$, $\beta\in \RR^{2}$}.
\ee


\section[Spectral stability conditions]{Bloch decomposition and 
spectral stability conditions}\label{bloch}
In co-moving coordinates, the linearized equation about $\bar u$ reads
\be\label{e:lin}
v_t=Lv:= \big((c-a)v\big)_x -v_{xxxx}-\epsilon v_{xxx} -\delta v_{xx},
\qquad a:= df(\bar u),
\ee
and the eigenvalue equation as
$Lv:=
- v_{xxxx}-\epsilon v_{xxx} -\delta v_{xx} + \big((c- av \big)_x =\lambda v.
$
Following \cite{G}, we define the one-parameter family of Bloch operators
\be \label{e:Lxi}
L_{\xi} := e^{-i \xi x} L e^{i \xi x},\quad\xi\in[-\pi,\pi)
\ee
operating on the class of $L^2$ periodic functions on $[0,X]$;
the $(L^2)$ spectrum
of $L$ is equal to the union of the
spectra of all $L_{\xi}$ with $\xi$ real with associated
eigenfunctions
$
w(x, \xi,\lambda) := e^{i \xi x} q(x, \xi, \lambda),
\label{e:efunction}
$
where $q$, periodic, is an eigenfunction of $L_{\xi}$.
By standard considerations,
the spectra of $L_{\xi}$
consist of the union of countably many continuous
surfaces $\lambda_j(\xi)$; see, e.g., \cite{G}.

Without loss of generality taking $X=1$,
recall now the {\it Bloch representation}
\be\label{Bloch}
u(x)=
\Big(\frac{1}{2\pi }\Big) \int_{-\pi}^{\pi}
e^{i\xi\cdot x}\hat u(\xi, x) d\xi
\ee
of an $L^2$ function $u$, where
$\hat u(\xi, x):=\sum_k e^{2\pi ikx}\hat u(\xi+ 2\pi k)$
are periodic functions of period $X=1$, $\hat u(\cdot)$
denoting with slight abuse of notation the Fourier transform of $u$
in $x$.
By Parseval's identity, the Bloch transform
$u(x)\to \hat u(\xi, x)$ is an isometry in $L^2$:
$
\|u\|_{L^2(x)}=\int_{-\pi}^\pi\int_0^1\left|\hat{u}(\xi,x)\right|^2dx~d\xi=
\|\hat u\|_{L^2(\xi; L^2(x))},
$
where $L^2(x)$ is taken on $[0,1]$ and $L^2(\xi)$ on $[-\pi,\pi]$.
More generally, for $q\le 2\le p$, $\frac{1}{p}+\frac{1}{q}=1$,
there holds the generalized Hausdorff--Young's inequality \cite{JZ}
\be\label{HY}
\|u\|_{L^p(x)}\le \|\hat u\|_{L^q(\xi;L^p(x))}.
\ee
The Bloch transform diagonalizes the periodic-coefficient operator $L$,
yielding the {\it inverse Bloch transform representation}
\be\label{IBFT}
e^{Lt}u_0=
\Big(\frac{1}{2\pi }\Big) \int_{-\pi}^{\pi}
e^{i\xi \cdot x}e^{L_\xi t}\hat u_0(\xi, x) d\xi \, .
\ee

Following \cite{JZ}, we assume along with (H1)--(H2) the
{\it strong spectral stability} conditions:

(D1) $\sigma(L_\xi) \subset \{ \hbox{\rm Re} \lambda <0 \} $ for $\xi\ne 0$.

(D2) $\hbox{\rm Re} \sigma(L_{\xi}) \le -\theta |\xi|^2$, $\theta>0$,
for $\xi\in \R$ and $|\xi|$ sufficiently small.

(\DDD) $\lambda=0$ is an eigenvalue
of $L_{0}$ of multiplicity $2$.\footnote{
The zero eigenspace of $L_0$,
corresponding to variations along the $3$-dimensional manifold
of periodic solutions in directions for which period does
not change \cite{JZ}, is at least $2$-dimensional
by (H2).
}

Assumptions (H1)-(H2) and (D1)--(\DDD) imply
\cite{JZ,BJNRZ3} that there exist $2$ smooth eigenvalues
\be\label{e:surfaces}
\lambda_j(\xi)= -i a_j \xi +o(|\xi|), \quad j=1,2
\ee
of $L_\xi$ bifurcating from $\lambda=0$ at $\xi=0$.
%
Following \cite{JZ}, we make the further nondegeneracy hypotheses:

(H3) The coefficients $ a_j$ in \eqref{e:surfaces} are distinct.

(H4) The eigenvalue $0$ of $L_0$ is nonsemisimple, i.e., $\dim
\ker L_0=1$.

\section[Whitham averaged equations]{Spectral
stability and the Whitham averaged equations}\label{s:whitham}

As noted in \cite{Se,JZ}, coefficients $a_j$ in \eqref{e:surfaces}
are characteristics of the $2\times 2$ Whitham averaged system
\ba\label{e:whitham}
M_t+ F_x&=0,\\
\omega_t + (c\omega)_x&=0
\ea
formally governing large-time
($\sim$ small frequency) behavior,
evaluated at the values $c,\omega$ of the background wave $\bar u$,
where $M(c,\omega)$ is the mean of $u$ over one period of
the periodic wave with speed $c$ and frequency $\omega=1/X$,
and $F(c,\omega)$ is the mean of $f(u)$.
Here, $\omega\sim \psi_x$, $c\sim -\psi_t/\psi_x$,
where $\psi$ denotes phase in the modulation approximation
\be\label{mod}
u(x,t)\approx \bar u(\psi(x,t)).
\ee

In the context of \eqref{e:KS}, thanks to the Galillean
invariance $x\to x-ct$, $u\to u+c$, \eqref{e:whitham}
reduces to
\ba\label{ksw}
c_t+ (H(\omega)-m(\omega)c)_x&=0,\\
\omega_t + (c\omega)_x&=0,
\ea
where $m(\omega)$ denotes the mean over one period of $u$
for a zero-speed wave of frequency $\omega$, and
$H(\omega)$ the mean of $u^2/2$, and in the classical situation
$\eps=0$ considered in \cite{FST}, to
$c_t+ (H(\omega))_x=0$,
$\omega_t + (c\omega)_x=0$,
which linearized about background values $c=0$, $\omega=\omega_0$, 
yields a {\it wave equation}
\be\label{e:wave}
\psi_{tt}+ \omega_0 dH(\omega_0)\psi_{xx}=0
\ee
so long as $dH(\omega_0)<0$.
Indeed, by odd symmetry, we may conclude in this case that the second-order
corrections $b_j$ in the further expansion
\be\label{further}
\lambda_j(\xi)=ia_j\xi- b_j \xi^2\cdots
\ee
of \eqref{e:surfaces} are equal,
hence $\lambda_j(\xi)$ agree to second order with the dispersion relations
of the {\it viscoelastic wave equation}
\be\label{e:dampedwave}
\psi_{tt}+ \omega_0 dH(\omega_0)\psi_{xx}=d(\omega_0)\psi_{txx},
\quad
d=2b_1=2b_2.
\ee

This recovers the formal prediction of ``viscoelastic behavior'' of
modulated waves carried out in \cite{FST} and elsewhere, or
``bouncing'' behavior of individual periodic cells. 
Put more concretely, \eqref{e:dampedwave} predicts that the
maxima of a perturbed periodic solution should behave approximately
like point masses connected by viscoelastic springs.
However, we emphasize that {\it such qualitative behavior}- in particular,
the fact that the modulation equation is of second order-  
{\it does not derive only
from Gallilean or other invariance of the underlying equations},
as might be suggested by early literature on the subject, 
{\it but rather from the more general structure of conservative} (i.e.,
divergence) {\it form} \cite{Se,JZ}.\footnote{As discussed further in \cite{Z},
conservation of mass lies outside the usual Noetherian formulation.}
Indeed, for any choice of $f$, $\lambda_j(\xi)$ may be seen to
agree to second order with the dispersion relation for an appropriate
diffusive correction of \eqref{e:whitham}, a generalized
viscoelastic wave equation.  See \cite{NR1,NR2} for further discussion
of 
Whitham averaged equations and their derivation.

\section{Linear estimates}

The main difficulty in obtaining linear estimates is that,
by (D3) and (H4), the zero eigenspace of $L_0$ has an associated
$2\times 2$ Jordan block.
This means that $e^{L_0 t}$ is not only neutral but grows as $O(t)$.
Viewed from the Bloch perspective, it means that the eigenprojections
of $L_\xi$ blow up as $\xi\to 0$.
Performing a careful spectral perturbation analysis, separating
out the singular part of the eigendecomposition of $e^{L_\xi t}$
in \eqref{IBFT},  and applying \eqref{HY},
as in Lemma 2.1, Prop. 3.3, and Prop. 3.4 of \cite{JZ},
we obtain the following detailed description of linearized behavior.

\begin{prop}\label{p:lin}
Under assumptions (H1)--(H4) and (D1)--(D3), the Green function $G(x,t;y)$ of \eqref{e:lin} decomposes
as
\be\label{dec}
G(x,t;y)=\bar{u}'(x)e(x,t;y)+\widetilde{G}(x,t;y),
\ee
where, for some $C>0$ and all $t>0$, $1\leq q\leq 2\leq p\leq\infty$ and $1\leq r\leq 4$,
\begin{align}
\left\|\int_{-\infty}^\infty\widetilde{G}(\cdot,t;y)f(y)dy\right\|_{L^p(\RM)}&\leq Ct^{-\frac{1}{4}\left(\frac{1}{2}-\frac{1}{p}\right)}
                          \left(1+t\right)^{-\frac{1}{4}\left(\frac{2}{q}-\frac{1}{2}-\frac{1}{p}\right)}\|f\|_{L^q\cap L^2}\label{finalGbd1}\\
\left\|\int_{-\infty}^\infty\partial_y^r\widetilde{G}(\cdot,t;y)f(y)dy\right\|_{L^p(\RM)}&\leq
         Ct^{-\frac{1}{4}\left(\frac{1}{2}-\frac{1}{p}\right)-\frac{r}{4}}
              \left(1+t\right)^{-\frac{1}{4}\left(\frac{2}{q}-\frac{1}{2}-\frac{1}{p}\right)-\frac{1}{2}+\frac{r}{4}}\|f\|_{L^q\cap L^2}\label{finalGbdyder}\\
\left\|\int_{-\infty}^\infty\partial_t\widetilde{G}(\cdot,t;y)f(y)dy\right\|_{L^p(\RM)}&\leq
        Ct^{-\frac{1}{4}\left(\frac{1}{2}-\frac{1}{p}\right)-1}
              \left(1+t\right)^{-\frac{1}{4}\left(\frac{2}{q}-\frac{1}{2}-\frac{1}{p}\right)+\frac{1}{2}}\|f\|_{L^q\cap L^2}\label{finalGbdtder}\, ,
\end{align}
$e(x,t;y)\equiv 0$ for $0\leq t\leq 1$, and
for all $t>0$, $1\leq q\leq 2\leq p\leq\infty$, $0\leq j,l,j+l\leq K$, and $1\leq r\leq 4$,
\begin{equation}\label{finalebds}
\begin{aligned}
\left\|\int_{-\infty}^\infty\partial_x^j\partial_t^l e(\cdot,t;y)f(y)dy\right\|_{L^p(\RM)}&\leq C\left(1+t\right)^{-\frac{1}{2}\left(\frac{1}{q}-\frac{1}{p}\right)
              -\frac{(j+k)}{2}}\|f\|_{L^q(\RM)}\\
\left\|\int_{-\infty}^\infty\partial_x^j\partial_t^l\partial_y^r e(\cdot,t;y)f(y)dy\right\|_{L^p(\RM)}&\leq C\left(1+t\right)^{\frac{1}{2}-\frac{1}{2}\left(\frac{1}{q}-\frac{1}{p}\right)
              -\frac{(j+k)}{2}}\|f\|_{L^q(\RM)} \, .
\end{aligned}
\end{equation}
Moreover,
for some constants $p_j$, and $a_j$ and $b_j$ as in \eqref{further},
\be\label{finest}
\Big\|e(\cdot,t;y)-
\sum_{j=1}^2 p_j  \errfn \left( \frac{\cdot-y-a_jt}{\sqrt{4b_jt}}, t \right)
\Big\|_{L^p}\le Ct^{-\frac{1}{2}(1-\frac{1}{p})} ,
\quad
t\ge 1.
\ee
\end{prop}

Defining $\psi:=-e$ and noting that $\bar u(x)+ \psi(x,t)\bar u'(x)\sim
\bar u(x+\psi(x,t))$, 
we see from \eqref{dec}--\eqref{finest} that linearized behavior indeed
agrees to lowest order with modulation by a phase function $\psi$
satisfying a generalized
viscoelastic wave equation obtained by diffusive correction 
of \eqref{e:whitham},
consisting of a {\it first-order hyperbolic--parabolic system}
in $\psi_x$, $\psi_t$.

This observation not only generalizes the second-order scalar
description \eqref{e:dampedwave} obtained in the special case 
$\eps=0$ \cite{FST}, but is in some sense more correct.
For, note that the second-order scalar description can be
a bit misleading as regards the assumption of initial data.
In particular, from the description \eqref{e:dampedwave},
one might be tempted to conclude that the linear response to a compactly
supported perturbation would consist of the D'Alembertian picture of
two approximately 
compactly supported wave forms in $\psi$ moving in opposite directions,
diffusing slowly at Gaussian rate.
Yet, the explicit bound \eqref{finest} shows that this is rather
a description of the derivative $\psi_x$!

Indeed, as described further in Section 1.2, \cite{JZN}, it is the
variables $\psi_x,\psi_t$ that are primary, rather than $\psi,\psi_t$ as
suggested by \eqref{e:dampedwave}, and it is these variables that
are related to the initial perturbation $(u-\bar u)|_{t=0}$.
{\it Thus, our analysis gives not only technical verification
of existing observations, but also new
intuition regarding the nature of modulated behavior.}


\section{Nonlinear stability}
Using the linear bounds of Prop. \ref{p:lin} together with 
nonlinear cancellation estimates as in \cite{JZ,JZN},
we obtain, finally, our main result describing nonlinear
behavior under localized perturbations.

\begin{theo}\label{main}
Assuming (H1)--(H4) and (D1)--(\DDD),
let $\bar u=(\bar \tau, \bar u)$ be a traveling-wave solution
of \eqref{e:gen}.
Then, for some $C>0$ and $\psi \in W^{2,\infty}(x,t)$,
\ba\label{eq:smallsest}
\|\tilde u-\bar u(\cdot -\psi-ct)\|_{L^p}(t)&\le
C(1+t)^{-\frac{1}{2}(1-1/p)}
\|\tilde u-\bar u\|_{L^1\cap H^K}|_{t=0},\\
\|\tilde u-\bar u(\cdot -\psi-ct)\|_{H^K}(t)&\le
C(1+t)^{-\frac{1}{4}}
\|\tilde u-\bar u\|_{L^1\cap H^K}|_{t=0},\\
\|(\psi_t,\psi_x)\|_{W^{K+1,p}}&\le
C(1+t)^{-\frac{1}{2}(1-1/p)}
\|\tilde u-\bar u\|_{L^1\cap H^K}|_{t=0},\\
\ea
and
\ba\label{eq:stab}
\|\tilde u-\bar u(\cdot-ct)\|_{ L^\infty}(t), \; \|\psi(t)\|_{L^\infty}&\le
C
\|\tilde u-\bar u\|_{L^1\cap H^K}|_{t=0}
\ea
for all $t\ge 0$, $p\ge 2$,
for solutions $\tilde u$ of \eqref{e:gen} with
$\|\tilde u-\bar u\|_{L^1\cap H^K}|_{t=0}$ sufficiently small.
In particular, $\bar u$ is nonlinearly bounded
$L^1\cap H^K \to L^\infty$ stable.
\end{theo}

Similarly as in the discussion of linear behavior,
we note that Theorem \ref{main} asserts
asymptotic $L^1\cap H^K\to L^p$
convergence of $\tilde u$ toward the modulated wave $\bar u(x-\psi(x,t))$,
but only bounded $L^1\cap H^K\to L^\infty$ stability about $\bar u(x)$,
a quite different picture from that suggested at first sight by 
\eqref{e:dampedwave}.

\section{Application to Kuramoto--Sivashinsky}\label{s:app}

We conclude our discussion by displaying some representative
traveling wave orbits and their associated spectrum, for
the case $\eps=0.2$, computed respectively using MATLAB
and the spectral Galerkin package SpectrUW.
These indicate, similarly as in the $\eps=0$ case studied
in \cite{CKTR,FST}, the existence of a band of spectrally
stable periodic traveling waves.
For related studies, and an animation of the spectral
evolution, see {\it Gallery}, http://www.math.indiana.edu/.
For more detailed numerical verification using the periodic
Evans function \cite{G}, see \cite{BJNRZ3}.

\begin{figure}[htbp]

\begin{center}
$
\begin{array}{ccc}
  \includegraphics[scale=.25]{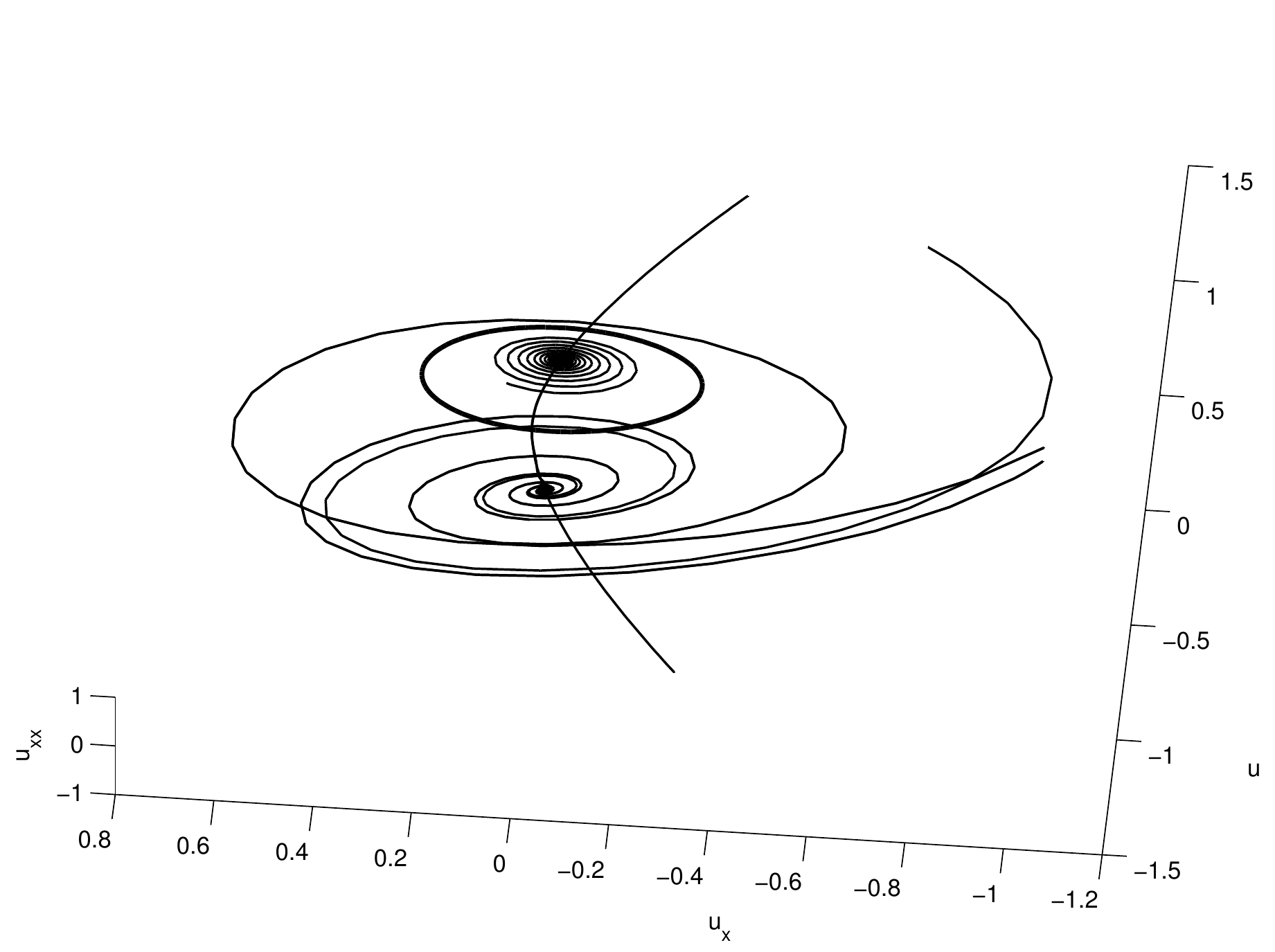}&\includegraphics[scale=0.25]{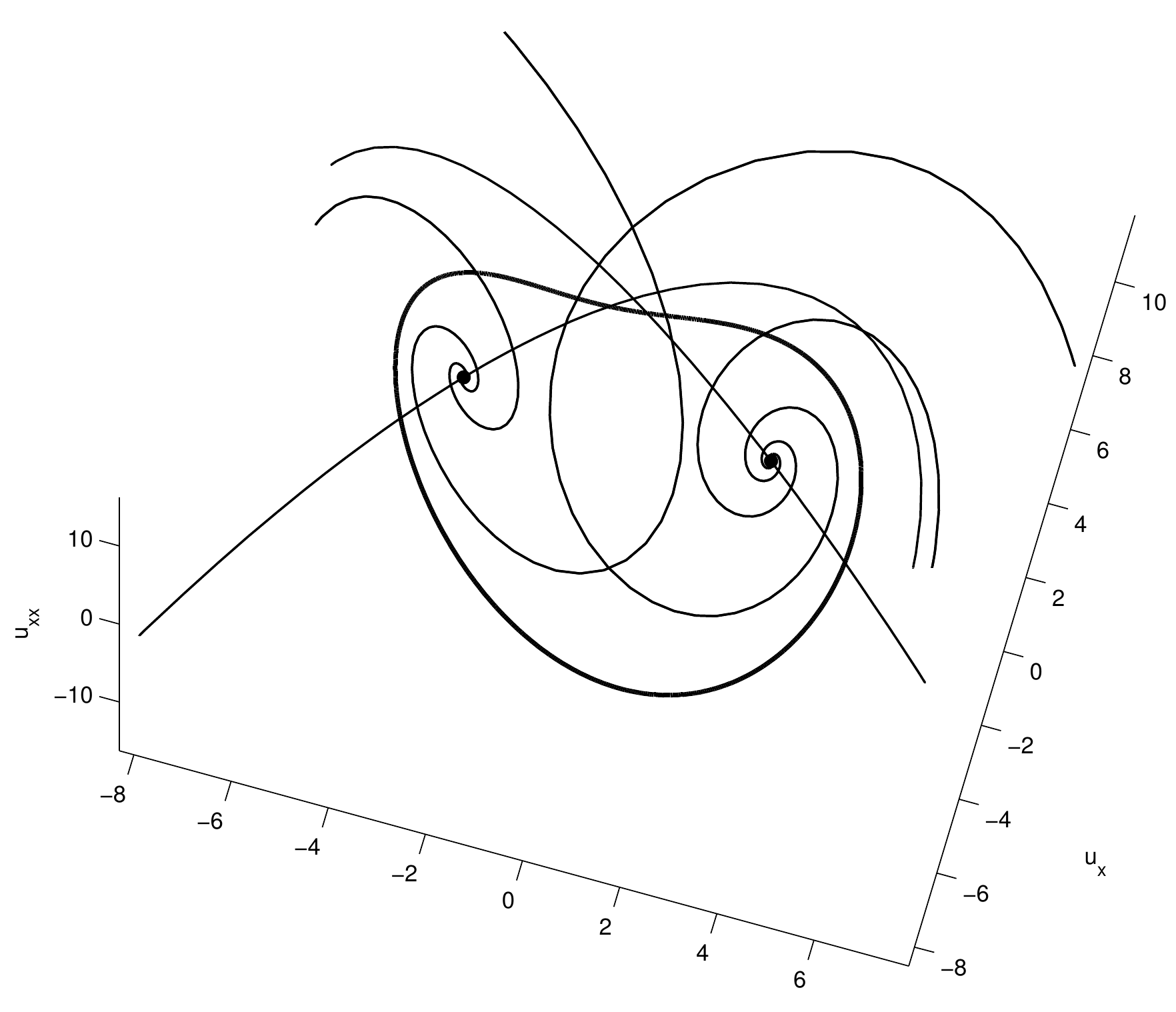}& \includegraphics[scale=.25]{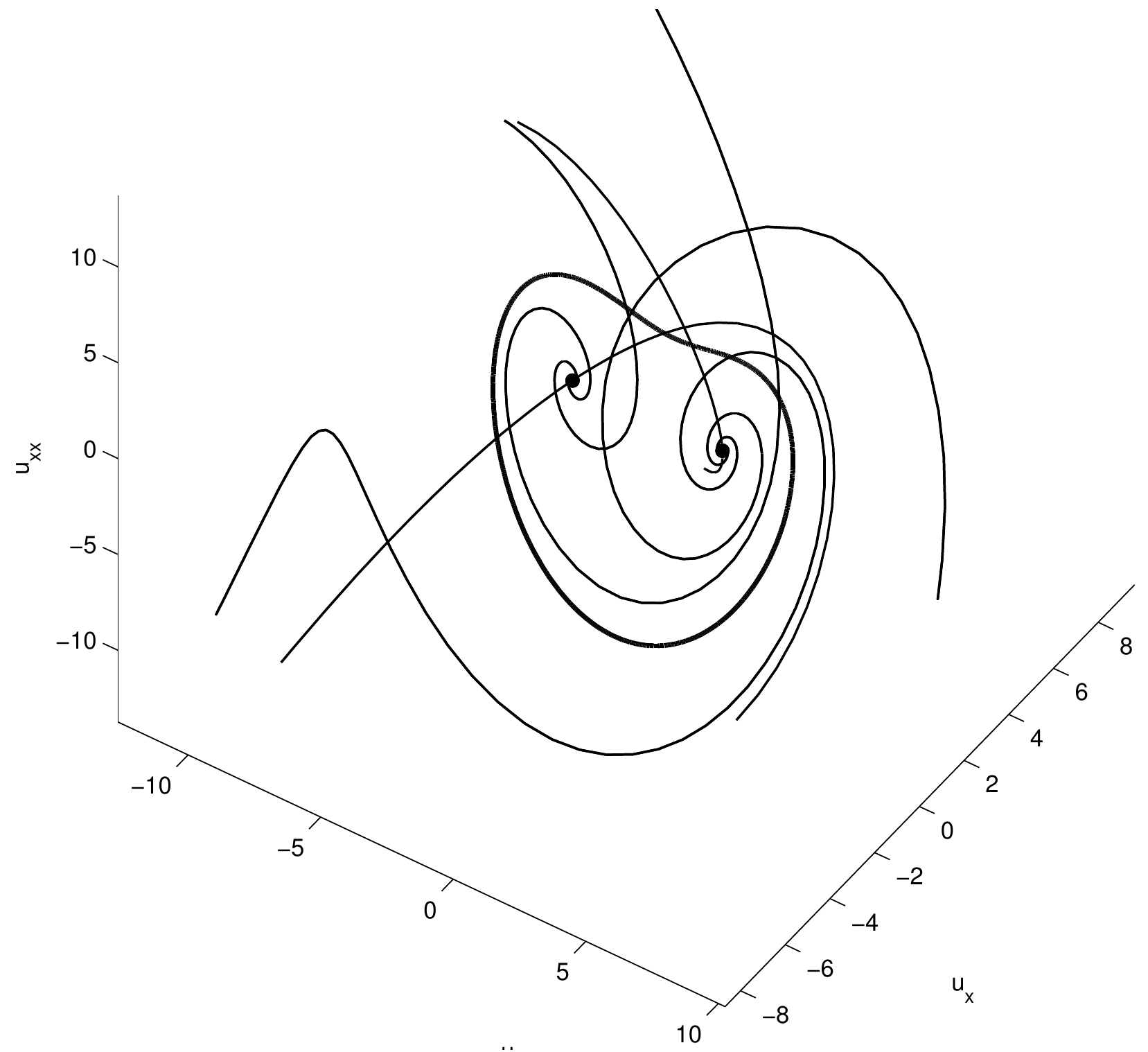}\\
  \includegraphics[scale=.25]{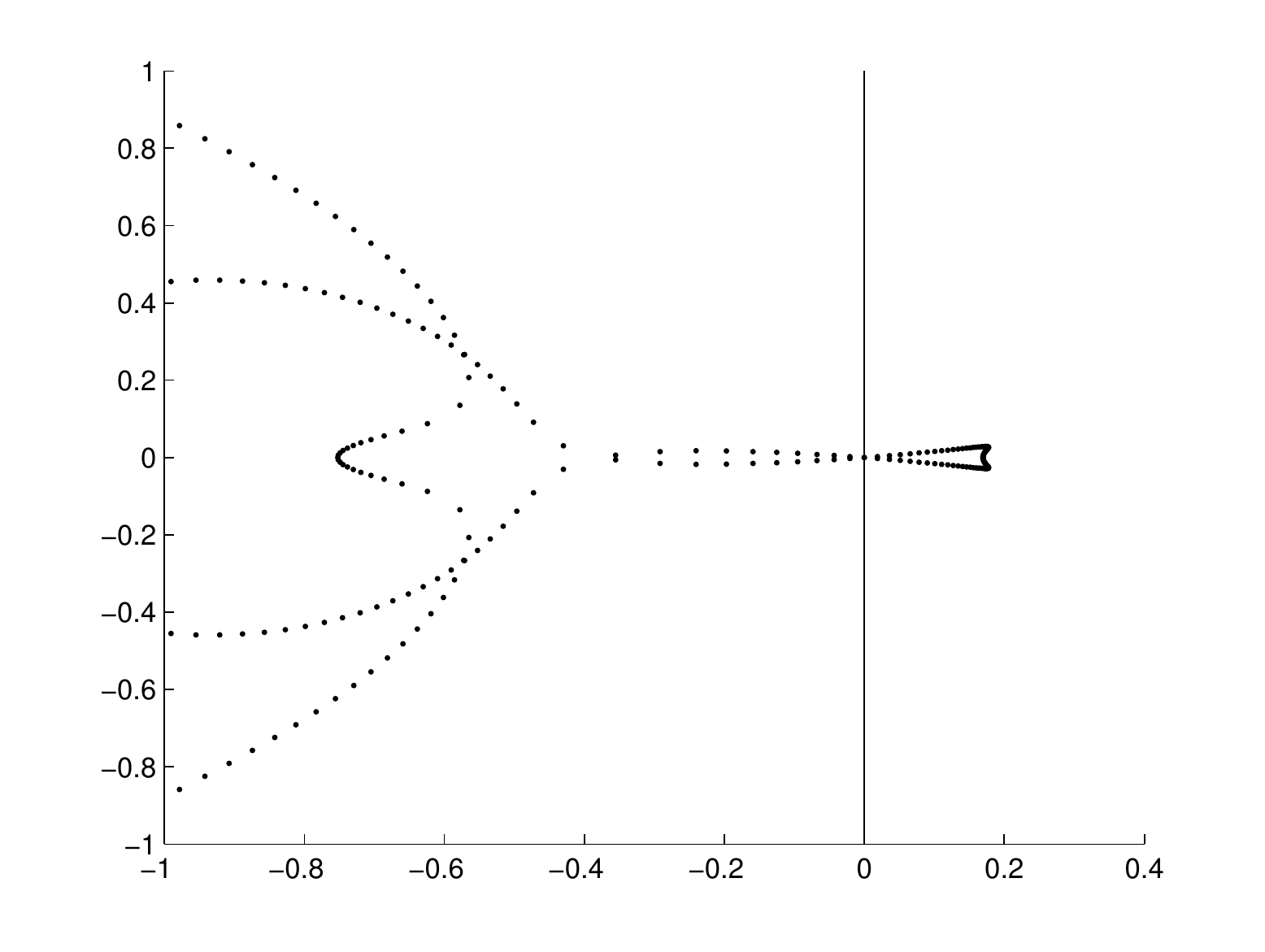}&\includegraphics[scale=0.25]{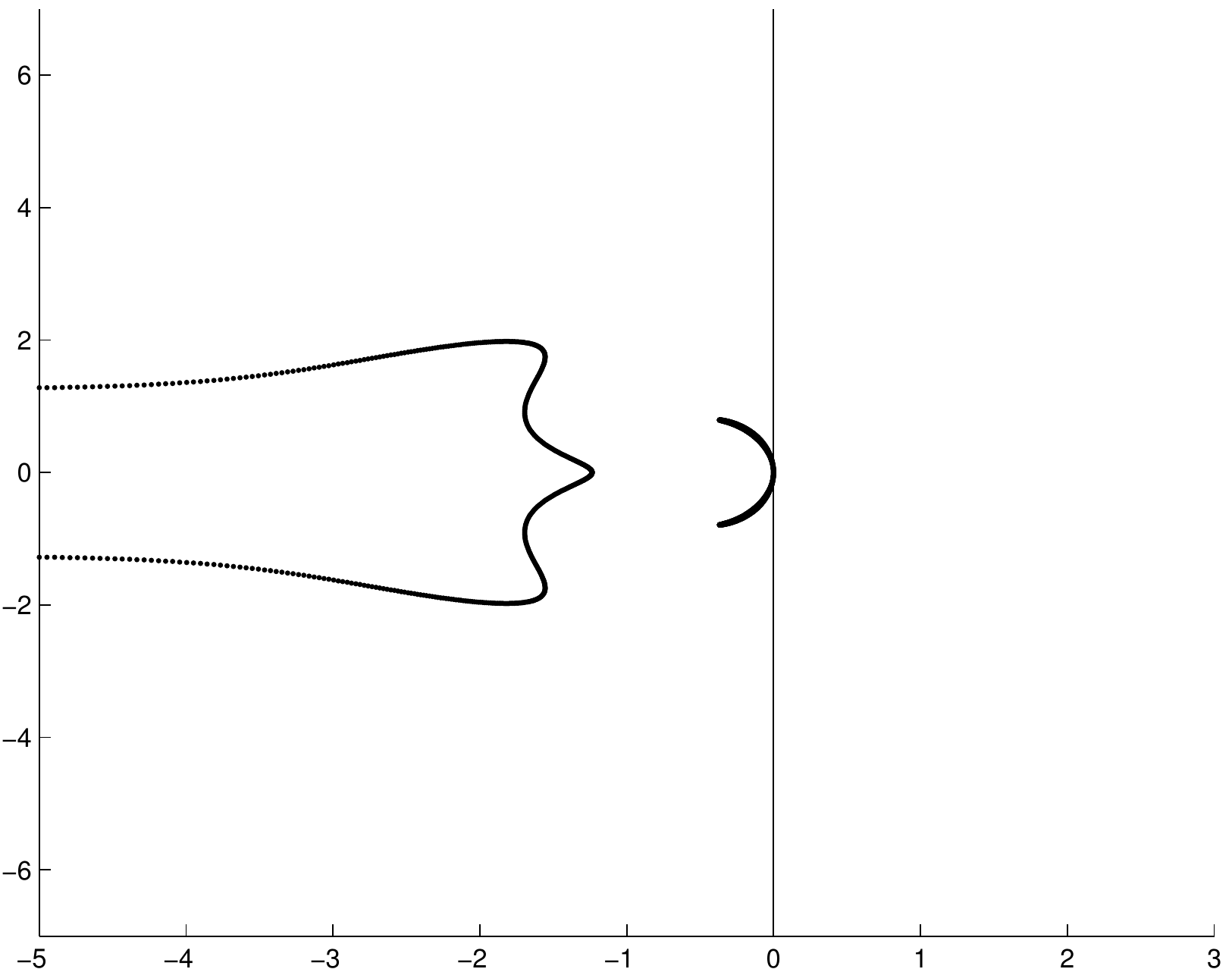}& \includegraphics[scale=.25]{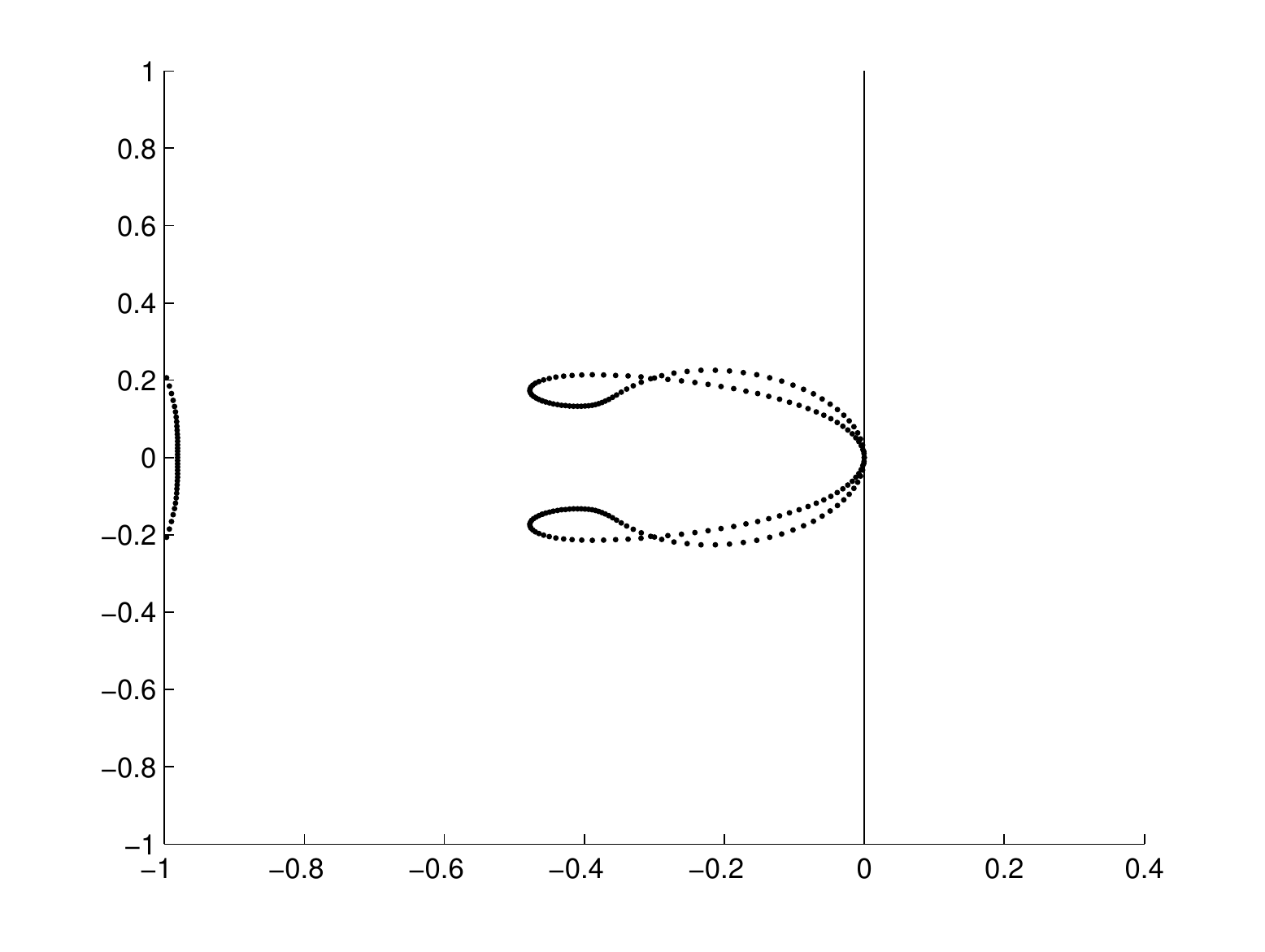}\\
\end{array}
$
\end{center}
\caption{Evolution of spectra as frequency $\omega=1/X$
varies with $c=0$ held constant.
Here $\eps=0.2$, $\delta=1$, $\gamma=1$.
Running left to right in the top row,
we see the evolution of the periodic orbits in the three-dimensional
phase portrait starting from Hopf
bifurcation and moving through the stable band of periodic waves,
with frames directly below depicting the spectrum of the corresponding
linearized operator about the wave.
Orbits were computed with MATLAB;
spectra were computed with the SpectrUW package
developed at University of Washington \cite{CDKK}.
}
\label{f:animate}
\end{figure}




\end{document}